\documentstyle[fullpage,amsfonts,12pt]{article}
\newcommand{\N}{{\Bbb N}}

\def\SBIMSMark#1#2#3{
 \font\SBF=cmss10 at 10 true pt
 \font\SBI=cmssi10 at 10 true pt
 \setbox0=\hbox{\SBF Stony Brook IMS Preprint \##1}
 \setbox2=\hbox to \wd0{\hfil \SBI #2}
 \setbox4=\hbox to \wd0{\hfil \SBI #3}
 \setbox6=\hbox to \wd0{\hss
             \vbox{\hsize=\wd0 \parskip=0pt \baselineskip=10 true pt
                   \copy0 \break%
                   \copy2 \break%
                   \copy4 \break}}
 \dimen0=\ht6   \advance\dimen0 by \vsize \advance\dimen0 by 8 true pt
                \advance\dimen0 by -\pagetotal
 \dimen2=\hsize \advance\dimen2 by .25 true in
     \setbox0=\hbox to 3.1 true in{
                \vbox to \ht6{\hsize=3 true in \parskip=0pt  \noindent  
 {\it Fund. Math.}~{\bf 155} (1998), 237--249.
                \vfill}}
  \ht0=0pt \dp0=0pt
 \ht6=0pt \dp6=0pt
 \setbox8=\vbox to \dimen0{\vfill \hbox to \dimen2{\copy0 \hss \copy6}}
 \ht8=0pt \dp8=0pt \wd8=0pt
 \copy8
 \message{*** Stony Brook IMS Preprint #1, #2 ***}
}

\begin{document}
 \title{Period Doubling, Entropy, and Renormalization}
 \author{Jun Hu$^{1}$ and
 Charles Tresser$^{2}$}
 \addtocounter{footnote}{1}\footnotetext
{Department of Mathematics,
 Graduate Center of CUNY, 33W 42nd Str., New York, NY 10036}
 \addtocounter{footnote}{1}\footnotetext{IBM, Po Box 218, Yorktown
 Heights, NY 10598}
\date{\relax}

 \maketitle

\SBIMSMark{1995/13}{October 1995}{}
\thispagestyle{empty}

 \begin{abstract}
 We show that in any family of
 stunted sawtooth maps, the set of maps whose set of periods is the
 set of all powers of $2$ has no interior point,
 {\em i.e.}, the combinatorial description of the boundary of chaos
 coincides with the topological description.
 We also show that, under mild assumptions, smooth multimodal maps
 whose set of periods is the set of all powers of $2$
 are infinitely renormalizable.
 \end{abstract}

\null
 \noindent {\bf 1 Introduction}

 \noindent
 The present work is motivated by the following
 folklore conjecture (see also [OT]):
 \begin{description}\item[Conjecture A.]
A real polynomial map of $f$ with set of periods (of its periodic orbits)
 \[P(f)=\{1, 2, 4, \cdots, 2^{i}, \cdots\}=\{2^{n}: n\in \N\} \]
 can be approximated by polynomials maps
with positive entropy and by polynomials maps with finitely
many periodic orbits.
 \end{description}

 \noindent
 This conjecture is now established for quadratic polynomials (as a
 consequence of [Su] or [La]) and work is in progress toward
 generalization for higher degree polynomials [Hu].
The interest in such a conjecture comes from Theorems A and B below
(see section 2.1) and the fact that {\em topological
entropy} (conceived as an invariant of topological conjugacy [AKM])
is also one way to measure the complexity of the dynamics
of a map (see section 2.1):
one is trying to describe how maps with simple dynamics can be
deformed to maps with complicated dynamics, or, as one says, chaotic
maps. Tradition,
as well as the availability in this framework of a greater set of
techniques, have put some emphasis on the particular case of polynomial
maps, as in Conjecture A. However, the problem of the transition to
chaos is more generally interesting in the category of smooth maps, in
particular smooth endomorphisms of the interval, for which we recall
the following:
 \begin{description}
 \item[Conjecture B.] All endomorphisms of the interval
 $f \in C^{k}(I),\, k\geq 1$,
 with $P(f)=\{2^{n}: n\in \N\}$ are on
 the boundary of chaos and on the boundary of the interior of the set of
 zero entropy.
 \end{description}

 \noindent
 We first show that any {\em stunted sawtooth map} (see section 2.2),
 whose set of periods is the
 set of all powers of $2$, can be approximated by
 stunted sawtooth maps with positive entropy and by stunted sawtooth maps
 with only finitely many periodic orbits (see section 2.3). This result
 solves the symbolic dynamic version of the above conjectures
in the sense that stunted sawtooth maps carry all possible
{\em kneading data} (see section 2.2) of multimodal maps.

 \noindent
 We also make a second step toward Conjecture A by
 proving that maps with $P(f)$ as above, which satisfy some smoothness
 conditions (and in particular polynomial maps), are {\em infinitely
 renormalizable} (see section 2.4).

 \noindent
Sections 3 and 4 contain proofs of the results formulated in section 2.

\noindent
{\em Acknowledgements.} Both authors would like to thank Dennis
Sullivan for his constant interest and encouragements. The
first author is grateful to the IMS at Stony Brook for support and
hospitality during his visit during the Spring of 1995. The second
author also benefited from regular visits to this institute.

 \vspace{.3in}
 \noindent {\bf 2 Preliminary definitions and results}

\bigskip
 \noindent {\bf 2.1 Topological entropy of one dimensional maps}

 \noindent
 A point $x$ is a $periodic$ $point$ of period $n$ of a map $f$ if
 $f^{i}(x)\neq x, 0<i<n$, and $f^{n}(x)=x$. The orbit of $x$ is then
 called a {\em periodic orbit} (of period $n$). If $n=1$, $x$ is
 called a $fixed$ $point$.

\noindent
 The {\em topological entropy} $h(f)$ of a continuous map $f$ on a
compact metric
 space $X$ with metric $d$ can be defined as follows [Bo].

\noindent
 Given $\epsilon >0$, $n\in Z_{+}$, we say a subset $S\subset X$ is
{\em $(n,\epsilon )$-separated} if
 \[x,y\in S, x\neq y\Rightarrow \exists m:0\leq m<n \;such\; that\;
 d(f^{m}(x),f^{m}(y))>\epsilon .\]
 Set $H(f,n,\epsilon )$ to be the maximal cardinality of $(n,\epsilon
 )$-separated
 sets, then
 \[h(f)=\lim _{\epsilon \rightarrow 0}\limsup _{n\rightarrow \infty }
 \frac{1}{n}\log H(f,n,\epsilon ).\]

\noindent
 For maps on an interval, the following result gives a necessary and
 sufficient condition for the positivity of topological entropy.
 \begin{description}
 \item[Theorem A.] ([BF], [Ml]) A continuous map of an interval
 to itself it has positive topological entropy if and only if it has
 a periodic point whose period is not a power of $2$.
 \end{description}

 \noindent
 {\bf Remark:} The ``if" part is from [BF], the ``only if"
part from [Ml].

\bigskip
 \noindent
 From Theorem A and [BlH], one gets the following.
 \begin{description}
 \item[Theorem B.]
 In the space $C^{k}(I), k\geq 1$, of $C^{k}$ endomorphisms
 of an interval $I$, if a map $f$ is on the boundary of positive
topological entropy then the set $P(f)$ of its
periods is $\{2^{n}:n\in \N\}$.
 The same is true for $f$ on the boundary of the interior of the set
 of maps with zero topological entropy.
 \end{description}

 \noindent
 {\bf Remark:} Conjectures A and B are about the converse of Theorem B.

\bigskip
 \noindent
We next give another necessary and sufficient condition for the
positivity of topological entropy, which will be an important tool
for us. This requires some more terminology.

 \noindent
 So let $f\in C^{0}(I)$ and $p$ be a fixed
point of $f$. A point $x$ of $I$ belongs to the {\em unstable manifold}
 $W^{u}(p,f)$
 of $p$ if, for every neighborhood of $V$ of $p$, $x\in f^{n}(V)$ for
 some positive integer $n$. It is easy to check that $W^{u}(p,f)$ is
 connected and invariant under $f$.
 A point $x\in I$ is a $homoclinic$ $point$ of $f$ if there is a periodic
 point $p$ of $f$ of period $n$ such that $x\neq p$, $x\in
 W^{u}(p,f^{n})$ and $f^{mn}(x)=p$ for some $m\in \N$ [Bl].

 \noindent
 \begin{description}
 \item[Theorem C]([Bl])
 A map $f\in C^{0}(I)$ has positive topological entropy if and only
 if it has a homoclinic point.
 \end{description}

 \noindent {\bf 2.2 Multimodal and stunted sawtooth maps}

 \noindent
Consider the continuous map $f:I\to I$ where $I=[c_0,c_{d+1}]$.
For $d\geq 0$, assume there are points
$c_i$, $0< i < d+1$ with $c_j< c_{j+1}$ for $0\leq j\leq {d}$
such that $f$ is monotone on each
{\em lap} $[c_j,c_{j+1}]$,
and not monotone on any segment of the form $[c_j,c_{j+2}]$. Such a map
is then called {\em d-modal} or {\em multimodal} with {\em modality}
$d$ (one says {\em amodal} if $d=0$, {\em unimodal} if $d=1$, and
then {\em bimodal}, and so on).
The maximal interval $[a_i,b_i]$ containing $c_i$ on which $f$
is constant is called
a {\em turning interval} and, more precilsely, a {\em plateau} if
$a_i<b_i$, and a {\em turning point} if $a_i=b_i$),

 \noindent
The {\em shape} of a $d$-modal map is the alternating sequence of $d+1$
signs, starting with either $+$ or $-$ according as the map is
increasing or decreasing on its initial lap.
By the {\em kneading data} associated with a $d$-modal map $f$ we will
mean its shape together with the collection of signs
$${\rm sgn}(f^n(c_i)-c_j)\in \{-1,0,1\}$$
for $n>0$ and $1\leq i,j\leq d$. The $i$-ordered collection of signs
with $j$ fixed is the $j^{\rm th}$ {\em kneading sequence} of $f$,
and the $j$-ordered collection of kneading sequences is the {\em kneading
invariant} of $f$ (for more on
kneading theory, we refer to [MiT] and [BORT]). One
might wish to first understand at the symbolic level some questions one
formulates for polynomials or smooth maps. It is
in fact more practical to consider continuous
families of $d$-modal maps rich enough
to exhibit all possible kneading data for $d$-modal maps, yet
significantly easier to study than smooth maps. Such families
exist, and we next recall the construction of one of them.

 \noindent
 By the {\em sawtooth map} of shape $s_{1}s_{2}\cdots s_{d+1},
 s_{i}\in \{+,-\}, s_{i+1}=-s_{i}$, $d\geq1$, we mean
 the unique map $S_{d}:I\rightarrow I$ which is piecewise linear with
 slope
 $s_{1}(d+1), s_{2}(d+1), \cdots s_{d+1}(d+1)$ ($d+1$ alternate values).
 This is a $d$-modal map with
 topological entropy $\log (d+1)$, the largest possible value
for $d$-modal maps.

\noindent
 Given
 any {\em critical value vector} $w=(w_{1}, w_{2}, \cdots , w_{d})$
satisfying
 $$(w_{j}-w_{j+1})\cdot s_{(j+1)}<0,\,w_{j}\in I,\,
j=1, 2, \cdots ,d\,,$$
we obtain the {\em stunted sawtooth map} $S_w$
 from $S_{d}$ by cutting off the tops and bottoms of the graph at heights
 $w_{j}, j=1,2,\cdots ,d$. The $d$-parameter family of stunted
sawtooth maps $S_w$ is complete in the following sense
 [MiT, DGMT].
 \begin{description}\item[Theorem D.] To any $d$-modal map $f$ there is a
 canonical $d$-modal stunted sawtooth map $S_w$ which has exactly the
 same kneading data as $f$.
 \end{description}

 \noindent {\bf 2.3 The first main result}
 \begin{description}\item[Theorem 1.]
 Suppose $S_{w}$ is a stunted sawtooth map with
 \[P(S_{w})=\{2^{i}:i\in \N\}.\]
 Then for any $\epsilon >0$, there exist $w^{'}$ and $w^{''}$ such that
 \[|w^{'}-w|<\epsilon ,\;\;\;\;\;|w^{''}-w|<\epsilon ,\]
 $h(S_{w^{'}})>0$ and $S_{w^{''}}$ has only finitely many periods.
 \end{description}
 \begin{description}\item[Corollary 1.]
 In the parameter space, the set $\{w: P(S_{w})=\{2^{n}: n\in \N\}\}$
 has no interior point, {\em i.e.}, the combinatorial and the topological
descriptions of the boundary of chaos coincide.
 \end{description}

 \noindent
{\bf Remark:} Let $f_i$, $i\in \{1,2,3\}$ be three
 $d$-modal maps of same shape with kneading invariants
$$K_i=\{K_{i,1},K_{i,2},\dots ,K_{i,d}\}\,.$$
Assume that $P(f_2)=\{2^{n}: n\in \N\}$ and that, with the usual order
on kneading sequences (see, {\em e.g.}, [BORT])
$$K_{1,j}< K_{2,j} < K_{3,j}\,\,{\rm if}\,\,s_j=+\,,$$
$$K_{1,j}> K_{2,j} > K_{3,j}\,\,{\rm if}\,\,s_j=-\,,$$
for $1\leq j\leq d$. It then also follows from Theorem 1
that $f_1$ has only finitely many periodic
orbits and that $f_3$ has positive topological entropy.

\bigskip
 \noindent
{\bf 2.4 Renormalization}

 \noindent
 Let $I$ be an interval. A map $f:I\rightarrow I$ is called
 $renormalizable$ if there exists a proper subinterval $J$ of $I$ and an
 integer $p$ such that

 (1) $f^{i}(J), i=0,1,\cdots ,p-1$, have no pairwise interior
 intersection,

 (2) $f^{p}(J)\subset J$.

 \noindent
 Then $f^{p}|_{J}:J\rightarrow J$ is called a $renormalization$ of $f$.
 A map $f: I\rightarrow I$ is {\em infinitely renormalizable} if there
 exist
 an infinite sequence $\{I_{n}\}_{n=1}^{\infty }$ of nested intervals
 and an infinite sequence $\{u(n)\}_{n=1}^{\infty }$ of integers such
 that $f^{u(n)}|_{I_{n}}: I_{n}\rightarrow I_{n}$ are renormalizations of
 $f$ and the length of $I_{n}$ tends to zero as $n\rightarrow \infty $.

 \noindent
 Another purpose of this paper is to prove that maps $f$
 with $P(f)=\{2^{i}:i\in \N\}$, which satisfy some smooth conditions,
 are infinitely renormalizable. To this end, we shall prove the
following abstract result for which we first recall some definitions.

 \noindent
 Let $f\in C^{0}(I)$. A periodic point $x$ of period $n$ is
  {\em attracting } (resp. {\em one-sided attracting
 }) if there exists a neighborhood (resp. a one-sided
 neighborhood) $U$
 of $x$ such that for any $y\in U$, $f^{nl}(x)\rightarrow x$ as
 $l\rightarrow \infty $. In both case we say the orbit of $x$ is
a $periodic$ $attractor$.
 An open interval $J\subset I$ is called a {\em wandering interval} of
 $f$ if

1) $f^{n}(J)\cap f^{m}(J) =\emptyset$ for any $n\neq m, n,m\in \N$
 and

 2) $f^{n}(J)$ does not converge to a periodic orbit.

 \begin{description}\item[Theorem 2]
 Assume the multimodal map $f: I\rightarrow I$ with $P(f)=\{2^{n}:
 n\in \N\}$ has no wandering intervals, no plateaus, and no more than
finitely many periodic attractors. Then $f$ is infinitely renormalizable.
 \end{description}

\noindent
Let $f: I\rightarrow I$ be a map which is not constant on any open set.
 We say that $f$ belongs to
 $\Gamma (2)$ if

 a) $f$ is $C^{2}$ away from the turning points;

 b) For every $x_{0}
 \in T_{f} $, there exists $\alpha >1 $, a neighborhood $U(x_{0}) \; of
 \; x_{0}  $
 and a $C^{2}$-diffeomorphism $\phi :U(x_{0}) \rightarrow
 (-1,1) $ such that $\phi
 (x_{0})=0$ and   \[f(x)=f(x_{0})\pm |\phi (x)|^{\alpha } ,\;\; \forall x
 \in  U(x_{0})	.\]

 \begin{description}\item[Theorem E.]$([Ma], [MMS])$

 (1) If $f\in \Gamma (2)$ then $f$ has no wandering
 intervals.

 (2) Any $f\in \Gamma (2)$ has at most finitely many periodic attractors.
 \end {description}

\noindent
Combining this result with Theorem 2 yields:

 \begin{description}\item[Theorem 3.]
Any $d$-modal map $f\in \Gamma (2)$, with $P(f)=\{2^{n}: n\in Z_{+}\}$
is infinitely renormalizable.
 \end{description}

 \noindent
In particular, we have also:

 \begin{description}\item[Theorem 4]
 Any real polynomial map $f$ with $P(f)=\{2^{n}:n\in \N\}$
 is infinitely renormalizable.
 \end{description}

\bigskip
\noindent
{\bf Remark:} Using results from [HS], the smoothness condition in (1) of
Theorem E can be relaxed, which allows a proof of
Theorem 3 with relaxed smoothness condition for unimodal maps.

\bigskip
\noindent
Using more language from kneading theory, one could formulate
a conjecture corresponding to Theorems 2 to 4 for general
renormalizations (not just those at the boundary of chaos).
Such a generalization completely escape the methods of the present paper.

\bigskip
 \noindent
{\bf 3 Proof of Theorem 1}
 \begin{description}\item[Lemma 1.]
 Let $f: I\rightarrow I$ be a multimodal map with
$P(f)=\{2^{i}:i\in Z_{+}\}$, and let $\Omega
 (f)$ be the set of accumulation points of the periodic
 points of $f$. Then no point in $\Omega (f)$ is periodic,
 so that $\Omega (f)$ is not a finite set.
\end{description}
 {\bf Proof:}
 Let $p\in \Omega (f)$ be a periodic point of period $2^{n}$. Denote
 $g=f^{2^{n}}$. Then $p$ is a fixed point of $g$. Look
at the map $g$ near
 $p$. Because $f$ hence $g$ has isolated turning intervals,
there are only three types of local behaviors for $g$.

 \noindent
 1. $g$ is monotone in a small neighborhood of $p$ (if $g$ is monotone
 reversing then $g^{2}$ is monotone preserving in a small neighborhood of
 $p$).

 \noindent
2. $p$ is in the interior or at the end of one of the plateaus of
 $g$.

 \noindent
3. $p$ is a turning point of $g$.

 \noindent
Look at three iterates of $g$
 or $g^{2}$ near $p$, one can easily see that neither
 $g$ nor $g^{2}$ has any periodic points with higher period in a
 small neighborhood of $p$. This contradicts that $p$ be an
 accumulation point of period-doubling periodic points of $f$. $\Box $

\bigskip
 \noindent
{\bf Remark:} Lemma 1 is false for continuous maps: examples with
$\Omega (f)$ reduced to a point are easily provided.

\bigskip
 \noindent
Let $\sum=\{0,1\}^{\N}$ and
let $\sigma$ stand for the {\em adding machine}, {\em i.e.},
the map $\sigma :\sum \rightarrow \sum $ defined by
$\sigma (x_i)_0^\infty = (y_i)_0^\infty $, where
$y_i=1-x_i$ if $x_j=1$ for all $j<i$ and $y_i=x_i$ otherwise.
 The following
result is proved in [M2].
 \begin{description}\item[Theorem F.]($[M2]$)
 Let $f\in C^{0}(I)$ be a continuous map with
 $P(f)=\{2^{n}: n\in \N^{+}\}$.
 Suppose that $K$ is an infinite closed invariant set of
$f$ which supports an ergodic $f$-invariant non-atomic probability
measure. Then there exists a continuous
map $h: K\rightarrow \sum $ such that
 $$h\circ f=\sigma \circ h.$$
 Furthermore $h^{-1}(s)$ contains at most two points for
 any point $s\in \sum $.
 \end{description}

\noindent
 In the proof of Theorem F (see [M2]), $K$ can be expressed as a disjoint
 union $K=K_{0}\cup K_{1}$,
 where the supporting intervals of $K_{0}$ and $K_{1}$ have disjoint
 interiors and
 $f(K_{i})=K_{1-i}$, where $i=0,1$. Hence $K_{0}$ and
 $K_{1}$ are invariant under $f^{2}$. $K_{0}$ and $K_{1}$ have
 the same bisections under $f^{2}$ and so on.
 Thus one can express $K$ as the disjoint unions
 \[K=\cup _{i=1}^{2^{n}}K^{(n)}_{i}\]
 for $n\in \N$. Therefore at least one point of each {\em fiber}
 $h^{-1}(s)$, $s\in \sum $, is recurrent.
 \begin{description}\item[Lemma 2.]
 Suppose that $S_{w}$ is a stunted sawtooth map with set of periods
 \[P(S_{w})=\{2^{i}:i\in Z_{+}\}.\]
 Let $\Omega (S_{w})$ be the set of accumulation points of
 period-doubling periodic
 points of $S_{w}$ and $\Sigma (S_{w})$ be the set of the endpoints of
 the interval $I$ and the closures of the turning intervals of $S_{w}$.
 Then the intersection of $\Omega (S_{w})$ and $\Sigma (S_{w})$ is not
 empty.
 \end{description}
 {\bf Proof:} Suppose that the intersection of $\Omega (S_{w})$ and
 $\Sigma (S_{w})$ is empty and let $\epsilon >0$ stand for the distance
$d(\Omega (S_{w}), \Sigma(S_{w}))$.
 Let $K\subset \Omega (S_{w})$ be as in Theorem F.
Clearly $d(K, \Sigma(S_{w}))>\epsilon$.
 Suppose that $x\in K$ is recurrent.
 There exists $l=2^{k}, k\geq 1$, (by Theorem F) such that
 $|f^{l}(x)-x|<\frac{1}{8}\epsilon $. Again by Theorem F, one can assume
 that $f^{l}(x)>x$. Let $V$ be the largest neighborhood of $x$ on which
 $f^{l}$ is monotone. The slope of $f^{l}$ on $V$ is $d^{l}$ and clearly
 $V\supset (x-\epsilon /d^{l}, x+\epsilon /d^{l})$.

\noindent
Assume first that $f^{l}$ is orientation
 preserving on $U$. Since $f^{l}(x)>x$ and $f^{l}(x-\frac{7}{8}\epsilon
 /d^{l})=f^{l}(x)-\frac{7}{8}\epsilon
 <x+\frac{1}{8}\epsilon -\frac{7}{8}\epsilon =x-\frac{3}{4}\epsilon
 <x-\epsilon
 /d^{l}$, there exists a point $p\in (x-\frac{7}{8}\epsilon /d^{l}, x)$
 such that $f^{l}(p)=p$. The unstable manifold $U(f^{l}, p)$ of
 $f^{l}$ at $p$ contains
 $(p, f^{l}(x)+\epsilon )$. We know that $f^{l}(f^{l}(x))<f^{l}(x)$
 by Theorem F. Let $W$ be the largest neighborhood of $f^{l}(x)$
 on which $f^{l}$ is monotone. Then $W\supset (f^{l}(x)-\epsilon /d^{l},
 f^{l}(x)+\epsilon /d^{l})$. No matter whether $f^{l}$ is preserving or
 reversing, $f^{l}(f^{l}(x)-\epsilon /d^{l})$ (or
 $f^{l}(f^{l}(x)+\epsilon /d^{l})$) is equal to $f^{l}(f^{l}(x))-\epsilon
 $, which is less than $f^{l}(x)-\epsilon <x+\frac{1}{8}\epsilon
 -\epsilon
 =x-\frac{7}{8}\epsilon <p$. This implies that there exists $y\in (x,
 f^{l}(x))\subset
 U(f^{l}, p)$ such that $y\neq p$ and $f^{l}(y)=p$, which is a homoclinic
 point, a contradiction. If $f^{l}$ is orientation
 reversing on $V$, one proceeds similarly.
 $\Box $

\bigskip
\noindent
 {\bf Proof of Theorem 1:} Since all turning points of $S_{w}$ (if
 any) are eventually fixed by $S_{w}$, by Lemma 1 they are not in
 $\Omega (S_{w})$. Let $\Lambda =E(S_{w})\bigcap
 \Omega (S_{w})$. For any $\epsilon >0$, we push the concave
 plateaus up a little bit and the convex plateaus down a little bit,
 if one of their end points belongs to	$\Lambda $, to get another
 stunted sawtooth map $S_{w^{'}}$ with $|w-w^{'}|<\epsilon $.
 If a periodic orbit of $S_{w}$ has no point in
 the interior of any plateau of $S_{w}$ then it is also a periodic orbit
 of $S_{w^{'}}$. Thus
 $S_{w}\subset S_{w^{'}}$so that
 the arguments in the proof of Lemma 2 show that $S_{w^{'}}$
 has a homoclinic point. Thus $S_{w^{'}}$ has positive
 topological entropy.

\noindent
 We next push all concave plateaus down a little bit and all
 convex plateaus up a little bit to get another stunted sawtooth map
 $S_{w^{''}}$ with $|w-w^{''}|<\epsilon $.
 Suppose that $P(S_{w^{''}})=\{2^{n}: n\in \N\}$.
 Then by the proof of Lemma 2, $S_{w}$ has positive
 topological entropy, a contradiction.
 $\Box $

 \vspace{.2in}
\goodbreak
 \noindent {\bf 4 Proof of Theorem 2}

 \begin{description}\item[Lemma 3.]
 Suppose $f\in C^{0}(I)$ has no attracting
 cycles and $P(f)=\{2^{n}: n\in \N\}$.
 Let $K$ be as in Theorem F, and $K=K_{0}\cup K_{1}$ with
 $f(K_{i})=K_{1-i}, i=0, 1$. Assume that $\sup \{x: x\in K_{0}\}<\inf
 \{x: x\in K_{1}\}$. Let $[K_{i}]$
 denote the smallest closed interval containing $K_{i}, i=0, 1$.
 Then for $i=0, 1$ there exists a periodic
 point $q$ of $f$ with period $n=1$ or $2$ in the gap between
 $[K_{0}]$ and $[K_{1}]$ such that $U(f^{2}, q)$ contains $[K_{i}]$.
 \end{description}

\noindent
 {\bf Proof:} By Theorem F, there exists a fixed point $p$ of $f$ in
 the gap between $[K_{0}]$ and $[K_{1}]$.
 For $i=1$, let $s=\inf \{x: x\in K_{1}\}$ and consider the map
 $f^{2}$ on the interval $[p,
 s]$. $p$ is fixed by $f^{2}$. Let $q$ be the largest
 fixed point of $f^{2}$ in $[p, s]$. Because of Lemma
 1, $q\neq s$ and then $f^{2}(y)>y$ for any $y\in (q, s)$. Clearly
 $f^{2}(K_{1})\subset K_{1}$ and $f^{2}(s)>s$.
 Therefore $s$, hence $[K_{1}]$by connectivity,
 is in $U(f^{2}, q)$.The case $i=0$ is treated similarly.
$\Box $.

\bigskip
\noindent
 From Theorem F and Lemma 3, one has
 \begin{description}\item[Corollary 2.]
 Suppose $f\in C^{0}(I)$ has no attracting
 cycles and $P(f)=\{2^{n}: n\in \N\}$.
 Let $K$ be as in Theorem F, with $K=\bigcup
 _{i=1}^{2^{n}}K_{i}^{(n)}$, $n\in \N$.
 Then for any $n\in \N$ there exists a periodic point $p$ of period
 $m$, where $m=2^{n}$ or $2^{n+1}$, whose orbit
is contained in the set $\bigcup
 _{i=1}^{2^{n}}[K_{i}^{(n)}]\setminus
 \bigcup _{i=1}^{2^{n+1}}[K_{i}^{(n+1)}]$ and such that
 $U(f^{m}, p)$ contains some $K_{i}^{(n+1)}$, where $1\leq i\leq
 2^{n+1}$.
 \end{description}
 \begin{description}\item[Lemma 4]
 Assume a continuous map $f: I\rightarrow I$ with $P(f)=\{2^{n}:
 n\in Z_{+}\}$ has no wandering interval,
no plateau, only finitely many
 turning points and finitely many periodic attractors.
Let $K$ be as in Theorem F.
 Then the semi-conjugacy $h$ in Theorem F is actually a conjugacy.
 \end{description}
 {\bf Proof:} Suppose that $h$ is not a conjugacy. Then
by Theorem F there
 exists a point $s\in \sum $ such that $h^{-1}(s)=\{x, y\}$, $x\neq y$.
 We claim that $h^{-1}(\sigma ^{n}(s))$ eventually contains a single
 point when $n$ is large enough. Let $I_{n}$ denote the supporting
 interval of
 $ h^{-1}(\sigma ^{n}(s))$, $n\geq 0$. Since there are only finitely many
 turning points and $I_{n}$, $n\geq 0$, are pairwise disjoint, there
 exists $m>0$ such that $I_{n}$ contains no turning points for any $n\geq
 m$. If the claim is false, then $I_{m}$ is a wandering interval, a
 contradiction. Therefore one can assume that $f(x)=f(y)$, where $x, y\in
 K$ and $x\neq y$. We separate our considerations into three cases.

 Case 1. Suppose that $f$ is monotone preserving in some neighborhoods
 of $x$ and $y$ (the proof is similar when $f$ is monotone reversing
 near $x$ and $y$). By Corollary 2, there exists a periodic point
 $p$ of $f$ of period $n$, $n=2^{k}$, $k>0$, such that the unstable
 manifold $U(f^{n}, p)$ of $f^{n}$ at $p$ contains the interval $[x, y]$,
 $p$ is not in $[x, y]$ and $p$ is close enough to $x$ or $y$. By
 continuity of $f$, when $p$ is close enough to $x$ (or $q$), $f(p)$
is very close to $f(x)=f(y)$. By the intermediate value theorem,
there are at least two
 points $a, b\in (x, y)$ such that $f(a)=f(b)=f(p)$. Then
 $f^{n}(a)=f^{n}(p)=p$. Hence $a$ is a homoclinic point, a contradiction.

 Case 2. Suppose that $f$ is monotone preserving in a neighborhood of $x$
 and monotone reserving in a neighborhood of $y$
(the proof is similar when
 the situation is reversed). Denote $x_{t}=f^{t}(x)=f^{t}(y)$, $t\in \N$.
 Let $\{c_{i}\}_{i\in \alpha }$ denote the set of
turning points of $f$, where
 $\alpha $ is a finite set. Clearly there exists $t_{0}>0$ such that
 $(x_{t}, c_{i})\bigcap K\neq \emptyset $ for any $t\geq t_{0}$,
 denoted by $K_{(t, i)}$. In the interval $(x, y)$ we select $v$ and $w$
 near $x$ and $y$ respectively with $f(v)=f(w)$. Because of the
 nonexistence of wandering intervals, the itineraries of $f(p)=f(q)$ and
 $f(v)=f(w)$ under $f$ are eventually different. Let $v_{t}=f^{t}(v)$,
 $t\in \N$. Then there exists $t>t_{0}$ such that there is at least one
 $c_{i}\in (x_{t}, v_{t})$. Clearly $K_{(t, i)}\subset
 (x_{t}, v_{t})$ and $f^{t}((x, v))=f^{t}((w, y))\supset (x_{t}, v_{t})$.
 From Theorem F and Corollary 2, there exists a periodic point $p$ of
 period $n=2^{k}$, $k>0$, where $p$ is not in $[x, y]$, such that $p$ is
 very
 near $x$ (or $y$), $f^{j}(p)$ is in the supporting interval $[K_{(t,
 i)}]$ and the unstable manifold $U(f^{n}, p)\supset [x, y]$.
 Clearly there exists a point $u\in (x, v)$ (or $u\in (w, y)$) such that
 $f^{t}(u)=f^{j}(p)$. Hence $u\in U(f^{n}, p)$, $u\neq p$ and
 $f^{nt}(u)=f^{n}(f^{j}(p))=p$. So $u$ is a homoclinic point, a
 contradiction.

 Case 3. If at least one of $x$ and $y$ is a turning point, the proof
 is just a slight modification of that of Case 1 or Case 2.   $\Box $

\bigskip
\noindent
 {\bf Proof of Theorem 2:} Let $K$ be as in Theorem F.
By Lemma 4, $f:
 K\rightarrow K$ is conjugate to $\sigma : \sum \rightarrow \sum $ by
 $h$. Separate the turning points of $f$ into two parts. Let $S_{1}$
 (resp. $S_{2}$) denote the turning points of $f$ contained (resp. not
 contained) in $K$. Denote $K=\cup _{i=1}^{2^{n}}K_{i}^{(n)}, n\in \N$.
 There exists $n_{0}\in \N$ such that for any $n>n_{0}$, $S_{2}\cap
 \cup _{i}^{2^{n}}[K_{i}^{(n)}]=\emptyset $. It is sufficient to prove
 for any $x\in S_{2}$ there exists $n_{0}$ such that for any $n>n_{0}$,
 $x$ is not in $\cup _{i}^{2^{n}}[K_{i}^{(n)}]$. Suppose
the contrary. Then there
 exists an infinite sequence of nested intervals $[K_{i(n)}^{(n)}]$
 containing $x$. Hence $x\in \cap _{n\in \N}[K_{i(n)}^{(n)}]$. Since the
 end points of $\cap _{n\in \N}[K_{i(n)}^{(n)}]$ are in the same fiber of
 $h$, which has to be a point, $x$ is in $K$, a contradiction.
 Now let $n>n_{0}$. Then $f$ maps each $[K_{i}^{(n)}]$ onto another
 $[K_{j}^{(n)}]$. Therefore $f$ is infinitely
 renormalizable. $\Box $
 
 \end{document}